\documentclass[12pt]{amsart}
\usepackage{amsmath}
\usepackage{amssymb}
\usepackage{bm}
\usepackage{graphicx}
\newtheorem{Theorem}{Theorem}
\newtheorem{Lemma}{Lemma}

\newtheorem{Corollary}{Corollary}

\begin{document}

\title[Graph Eigenfunctions and QUE]{Graph Eigenfunctions and Quantum Unique Ergodicity}
\author{Shimon Brooks and Elon Lindenstrauss}
\thanks{E.L. was
supported in part by NSF grants DMS-0554345 and DMS-0800345 and the Israel Science Foundation.}

\maketitle

{\em Abstract:}  We apply the techniques of \cite{meElon} to study joint eigenfunctions of the Laplacian and one Hecke operator on compact congruence surfaces, and joint eigenfunctions of the two partial Laplacians on compact quotients of $\mathbb{H}\times\mathbb{H}$.  In both cases, we show that quantum limit measures of such sequences of eigenfunctions carry positive entropy on almost every ergodic component.  Together with the work of \cite{Lin}, this implies Quantum Unique Ergodicity for such functions.
\medskip

\section{Introduction}
The purpose of this note is to point out an application of the techniques developed in the context of our work \cite{meElon} on eigenfunctions of large graphs to some Quantum Unique Ergodicity problems.

Our first result concerns certain compact hyperbolic surfaces $\Gamma\backslash\mathbb{H}$ of arithmetic congruence type.
One can consider more general $\Gamma$, but for concreteness and simplicity we restrict to the following situation.
Let $H$ be a quaternion division algebra over $\mathbb{Q}$, split over $\mathbb{R}$, and $R$ an order in $H$.  Fix an isomorphism $\Psi : H(\mathbb{R}) \cong \operatorname{Mat}_2(\mathbb{R})$.
For $\alpha\in R$ of positive norm $n(\alpha)$, we write $\underline{\alpha} = n(\alpha)^{-1/2}\Psi(\alpha)\in SL_2(\mathbb{R})$.  Set $\Gamma$ to be the image under $\Psi$ of the subgroup of norm $1$ elements of $R$.  As is well known, $\Gamma$ is discrete and co-compact in $SL_2(\mathbb{R})$, and so the quotient $X=\Gamma\backslash SL_2(\mathbb{R})$ is a $2$-to-$1$ cover of the unit cotangent bundle of a compact hyperbolic surface $M=\Gamma\backslash \mathbb{H}$.

Write $R(m)$ for the set of elements of $R$ of norm $m$, and define the Hecke operator\marginpar{}
$$T_m: f(x) \mapsto \frac{1}{\sqrt{m}}\sum_{\alpha\in R(1)\backslash R(m)} f(\underline{\alpha}x)$$
as the operator averaging over the Hecke points
$$T_m(x) = \{\underline{\alpha}x : \alpha\in R(1)\backslash R(m)\}$$
We will be interested in the case where $m=p^{k}$ are powers of a fixed prime $p$.  It is well known that $T_{p^k}$ is a polynomial in $T_p$; so in particular, eigenfunctions of $T_p$ are eigenfunctions of all $T_{p^k}$.  For all but finitely many primes, the points $T_{p^k}(x)$ form a $p+1$-regular tree as $k$ runs from $0$ to $\infty$; we will always assume that $p$ is such a prime.  We denote by $S_{p^k}$ the sphere of radius $k$ in this tree, given by Hecke points corresponding to the primitive elements of $R$ of norm $p^k$.

For any eigenfunction $\phi_j$ of the Laplacian $\Delta$ on $M$, normalized by $||\phi_j||_2=1$, one can construct a measure $\mu_j$ on $S^*M$ (which we view as a measure on the double cover $\Gamma \backslash SL_2(\mathbb{R})$) called the {\bf microlocal lift} of $\phi_j$ which is asymptotically invariant under the geodesic flow as the Laplace eigenvalue of $\phi_j$ tends to infinity. We shall use the variant of this construction used in \cite{LinHxH,Lin}, due to Wolpert \cite{Wolpert}, where $\mu_j=|\Phi_j|^2dVol$ for suitably chosen $\Phi_j \in L^2(S^*M)$ in the irreducible representation of $SL_2(\mathbb{R})$ on $\Gamma \backslash SL_2(\mathbb{R})$ generated by translates of $\phi_j$.
The construction satisfies that $\Phi_j$ is an eigenfunction of $T_p$ when $\phi_j$ is.
Since $\Delta$ commutes with all of the $T_p$, we may consider sequences $\{\phi_j\}$ of joint eigenfunctions, whereby each $\Phi_j$ is also an eigenfunction of $T_p$ (with the same eigenvalue as $\phi_j$).

\begin{Theorem}\label{main}
Let $p$ be a prime (outside the finite set of bad primes for $M$), and let $\{\phi_j\}_{j=1}^\infty$ be a sequence of $L^2$-normalized joint eigenfunctions of $\Delta$ and $T_p$ on $M$.  Then any weak-* limit point $\mu$ of the microlocal lifts $\mu_j$ has positive entropy on almost every ergodic component.
\end{Theorem}

Note that even without any Hecke operators, Anantharaman \cite{Anan} has shown (for general negatively curved compact manifolds) that any quantum limit has positive entropy, and this has been further sharpened in her joint work with Nonnenmacher \cite{AN} (see also \cite{AKN}). Hence the point of Theorem~\ref{main} is that it gives information on \emph{almost all ergodic components} of a quantum limit.
In view of the measure classification results of \cite{Lin}, this implies the following:

\begin{Corollary}
Let $\{\phi_j\}$ as above be a sequence of joint eigenfunctions of $\Delta$ and $T_p$.   Then the sequence $\mu_j$ converges weak-* to Liouville measure on $S^*M$.
\end{Corollary}

In \cite{Lin}, it was shown that if the sequence $\{\phi_j\}$ consists of joint eigenfunctions of $\Delta$ and {\em all} Hecke operators, then the $\mu_j$ converge to Liouville measure.  The assumption that $\phi_j$ is an eigenfunction of \underline{all} Hecke operators was used only to establish, through the work of Bourgain and the second named author \cite{BorLin}, that $\mu$ has positive entropy on a.e. ergodic component.
Hence Theorem~\ref{main} gives, at least formally, a strengthening of the result of \cite{Lin}.  It should be noted that if the multiplicities in the Laplace spectrum are uniformly bounded--- as is conjectured to be true in this setting--- then \cite{Lin} immediately implies QUE for any sequence of Laplace eigenfunctions, without making an additional assumption of Hecke invariance, and our result would be absorbed therein.  However, very little is known about these multiplicities at the present time.

Our methods also apply to the case of $M=\Gamma\backslash \mathbb{H}\times\mathbb{H}$, with $\Gamma$ a co-compact, irreducible lattice in $SL(2,\mathbb{R})\times SL(2,\mathbb{R})$.
Here we do not assume any arithmetic structure\footnote {By Margulis' Arithmeticity Theorem such lattices are necessarily arithmetic, though it is not known if they are necessarily of congruence type (which is necessary for the existence of Hecke Operators with good properties.)}; instead, we take the sequence $\{\phi_j\}$ to consist of joint eigenfunctions of the two partial Laplacians, each on the respective copy of $\mathbb{H}$ (equivalently, the $\phi_j$ are joint eigenfunctions of the full commutative algebra of $SL(2,\mathbb{R})\times SL(2,\mathbb{R})$-invariant differential operators).  The Laplacian on $M$ is the sum of the two partial Laplacians, and so the large eigenvalue limit for the Laplacian entails at least one of the two partial eigenvalues going to infinity (after passing to a subsequence, if necessary).

 Without loss of generality, we assume that the eigenvalues in the first coordinate are tending to infinity.
By usual semiclassical arguments (see \cite{LinHxH}), this means that the microlocal lift to $\Gamma\backslash SL(2,\mathbb{R})\times\mathbb{H}$ becomes invariant under the action of the diagonal subgroup $A$ of $SL(2,\mathbb{R})$ acting on the first coordinate.  By applying our analysis to the action of $SL(2,\mathbb{R})$ on the second coordinate--- which commutes with the action on the first coordinate--- ``in place" of the Hecke operator $T_p$, we are able to prove that any quantum limit of such a sequence must also carry positive entropy on a.e. ergodic component (with respect to the $A$-action on the first coordinate).  Thus, the result of \cite{Lin} again applies to show that $|\phi_j|^2dVol$ converges weak-* to the Riemannian measure $dVol$ on $M$.
The argument is analogous to the one presented here for the rank-one arithmetic case in Theorem~\ref{main}, and we intend to supply complete proofs of both results in a forthcoming paper.

The results of \cite{BorLin} were generalized by Silberman and Venkatesh \cite{Lior-Akshay, Lior-Akshay2} who, using the measure classification results of \cite{EKL}, were able to extend the QUE results of \cite{Lin} to quotients of more general  symmetric spaces (they also had to develop an appropriate microlocal lift). It is likely possible to extend the techniques of this paper to their context.

Regarding finite volume arithmetic surfaces such as $SL(2,\mathbb{Z})\backslash\mathbb{H}$, it was shown in \cite{Lin} that any quantum limit
has to be a scalar multiple of the Liouville measure, though not necessarily with the right scalar, and similar results can be provided by  our techniques using a single Hecke operator. Recently, Soundararajan \cite{Sound} has given an elegant argument that settles this escape of mass problem and (in view of the results of \cite{Lin}) shows that the only quantum limit is the normalized Liouville measure. An interesting open question is whether our $p$-adic wave equation techniques can be used to rule out escape of mass using a single Hecke operator. We also mention that Holowinsky and Soundararajan \cite{Holowinsky-Sound} have recently developed an alternative approach to establishing Arithmetic Quantum Unique Ergodicity for joint eigenfunctions of all Hecke operators. This approach \emph{requires} a cusp, and is only applicable in cases where the Ramanujan Conjecture holds; this conjecture is open for the Hecke-Maass forms, but has been established by Deligne for holomorphic cusp forms --- a case which our approach does not handle.

\section{The Propagation Lemma}

The following lemma proved along the lines of \cite{meElon} is central to our approach:

\begin{Lemma}\label{operator}
Let $\eta>0$.
For any sufficiently large $N\in \mathbb{N}$ (depending on $\eta$), and any $T_p$ eigenfunction $\Phi_j$,
there exists a convolution operator $K_N$ on $S^*M$ satisfying:
\begin{itemize}
\item{$K_N(\delta_x)$ is supported on the union of Hecke points $y\in T_{p^j}(x)$ up to distance $j\leq N$ in the Hecke tree, and is  constant on the spheres $S_{p^j}(x)$.}
\item{$K_N$ has matrix coefficients bounded by $O(p^{-N\delta})$, in the sense that for any $x\in S^*M$
$$|K_N(f)(x)|\lesssim p^{-N\delta} \sum_{j=0}^{N} \sum_{y\in S_{p^j}(x)} f(y)$$
}
\item{Any $T_p$ eigenfunction is also an eigenfunction of $K_N$, of eigenvalue $\geq -1$.  Moreover, $\Phi_j$ has $K_N$-eigenvalue $>\eta^{-1}$.}
\end{itemize}
\end{Lemma}

Lemma~\ref{operator} is based on the well known connection between Hecke operators and Chebyshev polynomials. A way to derive these which we have found helpful is via the following $p$-adic wave equation for functions on~$\mathcal{G}$
\begin{eqnarray*}
\Phi_{n+1} & = & \frac{1}{2}T_p\Phi_n - \left(1-\frac{T_p^2}{4}\right)\Psi_n\\
\Psi_{n+1} & = & \frac{1}{2}T_p\Psi_n + \Phi_n
\end{eqnarray*}
which is a discrete analog of the non-Euclidean wave equation (more precisely, of the unit time propagation map for the wave equation) on $\mathbb{H}$.  For initial data $(\Phi_0, \Psi_0)\in L^2(\mathcal{G})\times L^2(\mathcal{G})$, the solution to this equation is given by the sequence
\begin{eqnarray*}
\Phi_n & = & P_n\left[\frac{1}{2}T_p\right] \Phi_0 - \left(1-\frac{T_p^2}{4}\right) Q_{n-1}\left[\frac{1}{2}T_p\right]\Psi_0\\
\Psi_n & = & P_n\left[\frac{1}{2}T_p\right] \Psi_0 + Q_{n-1}\left[\frac{1}{2}T_p\right]\Phi_0
\end{eqnarray*}
where $P$ and $Q$ are Chebyshev polynomials of the first and second kinds, respectively, given by
\begin{eqnarray*}
P_n(\cos{\theta}) & = & \cos{n\theta}\\
Q_{n-1}(\cos{\theta}) & = & \frac{\sin{n\theta}}{\sin{\theta}}
\end{eqnarray*}
This can be proved directly by induction, using the recursive properties of the Chebyshev polynomials.

Suppose we take initial data $(\delta_0, 0)$.  The solution to the $p$-adic wave equation is then $\{(P_n[\frac{1}{2}T_p]\delta_0, Q_{n-1}[\frac{1}{2}T_p]\delta_0)\}$.  On the other hand, we can compute the explicit solution inductively; looking at the first coordinate, we get the following ``Propagation Lemma" on the tree:
\begin{Lemma}\label{propagation}
Let $\delta_0$ be the delta function at $0$ in the $p+1$-regular tree.  Then for $n$ even, we have
\begin{eqnarray*}
P_n\left[\frac{1}{2}T_p\right]\delta_0(x) & = & \left\{ \begin{array}{ccc} 0  & \quad & |x|  \text{ odd } \quad \text{or} \quad |x|>n\\ \frac{1-p}{2p^{n/2}} & \quad & |x|<n \quad \text{and} \quad |x| \text{ even } \\ \frac{1}{2p^{n/2}} & \quad & |x| = n \end{array}\right.
\end{eqnarray*}
In particular, we have
$$P_n\left[\frac{1}{2}T_p\right]\delta_0(x) \lesssim p^{-n/2}$$
\end{Lemma}

We now have a description of the $p$-adic wave propagation in both spectral and spacial terms, which we will use to construct our desired radial kernel $K_N$ on the Hecke tree.  As is already evident in Lemma~\ref{propagation}, it will be convenient to write the $T_p$-eigenvalues as $2\cos(\theta)$, where
\begin{itemize}
\item{The tempered spectrum is parametrized by $\theta\in [0, \pi]$.}
\item{The positive part of the untempered spectrum has $i\theta \in (0,\log{\sqrt{p}})$. }
\item{The negative part of the untempered spectrum has $i\theta+\pi \in (0,\log{\sqrt{p}})$. }
\end{itemize}

Consider first the case where $\Phi_j$ has $T_p$-eigenvalue $2$, or $\theta=0$.  Denoting the Fej\'er kernel of order $M$ by $F_M$, we set the spherical transform of $K_N$ to be
$$h_{K_N}(\theta) = F_M(q\theta)-1$$
Now $F_M(0)=M$ and $F_M$ is non-negative, so the third condition of Lemma~\ref{operator} is satisfied on the tempered spectrum, as long as $M>\eta^{-1}+1$.  Moreover, we can write
$$F_M(q\theta)-1 = \sum_{j=1}^M \frac{2(M-j)}{M} \cos{(jq\theta)}$$
and observing that $\cos(jq\theta)>\cos(0)$ on the entire untempered spectrum as long as $q$ is even, we see that the third condition holds on the full spectrum.
We also observe that
\begin{eqnarray*}
||K_N||_{L^1(\mathcal{G})\to L^\infty(\mathcal{G})}
& \leq & \sum_{j=1}^M 2\left|\left|P_{jq}\left[\frac{1}{2}T_p\right]\right|\right|_{L^1(\mathcal{G})\to L^\infty(\mathcal{G})}\\
& \lesssim & \sum_{j=1}^M p^{-jq/2}\\
& \lesssim & p^{-q/2}
\end{eqnarray*}
which satisfies the second condition of Lemma~\ref{operator}, as long as $q\geq 2N\delta$.  Moreover, since each $P_{jq}[\frac{1}{2}T_p]\delta_0$ vanishes outside the ball of radius $Mq$, the first condition is satisfied once $Mq\leq N$.  So we may take $M=\lceil\eta^{-1}\rceil+1$, and $q=2\lfloor N/2M\rfloor$, which yields $\delta = \lfloor q/2N\rfloor \gtrsim \eta$.  The same kernel also works for $\theta=\pi$ or untempered $\theta$.

We must now  consider $\theta\in (0,\pi)$.  By Dirichlet's Theorem, we may find a $q' < N\eta$ such that $q'\theta$ is close to $0$; so close, in fact, that we may take as large an even multiple $q$ of $q'$ as necessary to insure that, for a suitable constant $c$, we have
$$cN\eta^2 < q < 2N\eta$$
while $q\theta$ is still close enough to $0$ that we have
$F_{2M}(q\theta) > \eta^{-1}+1$.
Thus in this case we can take $\delta\gtrsim c\eta^2$, and Lemma~\ref{operator} follows.

\section{Proof of the Main Result}
The following estimate can be derived using the techniques of \cite{BorLin}, specifically Lemmas~3.1 and 3.3 there (much more general statements of this type by Silberman and Venkatesh can be found in \cite{Lior-Akshay2}).
\begin{Lemma}\label{disjoint}
For $\tau$ fixed but small enough, there exists a constant $c$ (depending only on $\tau$), such that for any $x,z\in X$, and any
$\epsilon<cp^{-2N}$, the tube $zB(\epsilon, \tau)$ contains at most $O(N)$ of the Hecke points
$$\{\underline{\alpha}x: \alpha\in T_{p^{j}}(x) \text{ for some } j\leq N  \}$$
\end{Lemma}

\noindent
{\em Sketch of proof of Theorem~\ref{main}:} Take any sufficiently fine partition $\mathcal{P}$ of $S^*M$, and consider its refinement under the time one geodesic flow.  Any partition element of the refinement is contained in a union of $O_c(1)$ tubes of the form $xB(cp^{-2N},\tau)$ for some $x\in S^*M$.  We assume that $c$ is sufficiently small, according to Lemma~\ref{disjoint}, independent of $\eta$ and $N$.

Take a collection $\{E_1, E_2, \ldots, E_K\}$ of distinct partition elements of the $\lfloor 2N\log{p}\rfloor$-th refinement of $\mathcal{P}$, of cardinality $K$, and set $\mathcal{E}=\bigcup_{k=1}^K E_k$ to be their union.  Let $1_{E_k}$ denote the characteristic function of each $E_k$, and similarly $1_{\mathcal{E}} = \sum_{E_k\subset\mathcal{E}} 1_{E_k}$.
To each $E_k$ we associate, as above, $O_c(1)$ tubes $B_{k,l} = x_{k,l}B(cp^{-2N},\tau)$ whose union contains $E_k$.

Now assume that $\mu(\mathcal{E})>\eta$; by definition, this implies that there exists a $j$ such that $\mu_j(\mathcal{E})=||\Phi_j1_\mathcal{E}||_2^2 > \eta$ as well.   Consider the correlation
$$\langle K_N (\Phi_j1_{\mathcal{E}}),  \Phi_j1_{\mathcal{E}}\rangle$$
We will estimate this in two different ways.  First, since $K_N$ has small matrix coefficients, we have
\begin{eqnarray*}
\langle K_N (\Phi_j1_{\mathcal{E}}),  \Phi_j1_{\mathcal{E}}\rangle_{L^2(S^*M)}
& = & \sum_{k=1}^K \langle K_N(\Phi_j1_\mathcal{E}), \Phi_j1_{E_k}\rangle_{L^2(E_k)}\\
& \leq & \sum_{k=1}^K ||K_N(\Phi_j1_\mathcal{E})||_{L^2(E_k)}||\Phi_j1_{E_k}||_2\\
& \leq & \max_{1\leq k\leq K}||K_N(\Phi_j1_\mathcal{E})||_{L^2(E_k)}\sum_{k=1}^K ||\Phi_j1_{E_k}||_2\\
& \leq & \sum_{i=1}^K \max_{1\leq k\leq K}||K_N(\Phi_j1_{E_i})||_{L^2(E_k)} \sum_{k=1}^K ||\Phi_j1_{E_k}||_2\\
& \lesssim & \sum_{i=1}^K \max_{1\leq k\leq K} N p^{-\delta N} ||\Phi_j1_{E_i}||_2 \sum_{k=1}^K ||\Phi_j1_{E_k}||_2
\end{eqnarray*}
by Lemmas~\ref{disjoint} and \ref{operator}, since each $E_i$ contributes at most $O(N)$ terms of size $p^{-\delta N}$ to $||K_N(\Phi_j1_{E_i})||_{L^2(E_k)}$---  were there $x\in E_k$ having $CN$ Hecke points in $\bigcup B_{i,l}$, then there would be have to be $\gtrsim CN$ Hecke points in a single $B_{i,l}$, all at distance $\leq N$ in the Hecke tree; this contradicts Lemma~\ref{disjoint} once $C$ is large enough.
Therefore,
\begin{eqnarray}
\langle K_N (\Phi_j1_{\mathcal{E}}),  \Phi_j1_{\mathcal{E}}\rangle
& \lesssim & \sum_{i=1}^K \max_{1\leq k\leq K} Np^{-\delta N} ||\Phi_j1_{E_i}||_2 \sum_{k=1}^K ||\Phi_j1_{E_k}||_2\nonumber\\
& \lesssim & Np^{-\delta N} \left(\sum_{k=1}^K ||\Phi_j1_{E_k}||_2 \right)^2\nonumber\\
& \lesssim & Np^{-\delta N} \sum_{k=1}^K||\Phi_j1_{E_k}||_2^2 \cdot K\nonumber\\
& \lesssim & Np^{-\delta N} K\label{pointwise}
\end{eqnarray}
since $\sum_{k=1}^K ||\Phi_j1_{E_k}||_2^2 \leq ||\Phi_j||_2^2 = 1$.

On the other hand, we can decompose $\Phi_j1_\mathcal{E}$ spectrally into an orthonormal basis (of $L^2(S^*M)$) of $T_p$ eigenfunctions $\{\psi_i\}\ni \Phi_j$, which {\em a fortiori} also diagonalize $K_N$, and notice that
\begin{eqnarray*}
\Phi_j1_\mathcal{E} & = & \langle \Phi_j1_\mathcal{E}, \Phi_j\rangle \Phi_j + \sum_{\psi_i\neq \Phi_j} \langle \Phi_j1_\mathcal{E}, \psi_i\rangle \psi_i\\
& = & ||\Phi_j1_\mathcal{E}||_2^2 \Phi_j + \sum_{\psi_i\neq \Phi_j} \langle \Phi_j1_\mathcal{E}, \psi_i\rangle\psi_i
\end{eqnarray*}
with
\begin{eqnarray*}
\sum_{\psi_i\neq \Phi_j} |\langle \Phi_j1_\mathcal{E}, \psi_i\rangle|^2 & = & ||\Phi_j1_\mathcal{E}||_2^2 - |\langle \Phi_j1_\mathcal{E}, \Phi_j\rangle|^2\\
& = & ||\Phi_j1_\mathcal{E}||_2^2 - ||\Phi_j1_\mathcal{E}||_2^4\\
& < & ||\Phi_j1_\mathcal{E}||_2^2 (1-\eta)
\end{eqnarray*}
by the assumption that $||\Phi_j1_\mathcal{E}||_2^2 > \eta$.

Now by Lemma~\ref{operator}, since $\{\psi_i\}$ diagonalizes $K_N$, and the $K_N$ eigenvalue of each $\psi_i$ is at least $-1$, we have
\begin{eqnarray}
\langle K_N (\Phi_j1_{\mathcal{E}}),  \Phi_j1_{\mathcal{E}}\rangle
& = & \sum_{\psi_i}|\langle \Phi_j1_\mathcal{E}, \psi_i\rangle|^2\langle K_N\psi_i, \psi_i\rangle\nonumber\\
& \geq & |\langle \Phi_j1_\mathcal{E}, \Phi_j\rangle|^2 \langle K_N\Phi_j, \Phi_j\rangle - \sum_{\psi_i\neq \Phi_j} |\langle \Phi_j1_\mathcal{E}, \psi_i\rangle|^2\nonumber\\
& > & ||\Phi_j1_\mathcal{E}||_2^4\langle K_N\Phi_j, \Phi_j\rangle - ||\Phi_j1_\mathcal{E}||_2^2 (1-\eta)\nonumber\\
& > & ||\Phi_j1_\mathcal{E}||_2^2 (||\Phi_j1_\mathcal{E}||_2^2\cdot \eta^{-1} - (1-\eta))\nonumber\\
& > & \eta(\eta\cdot\eta^{-1} -1 +\eta) = \eta^2 >0\label{spectral}
\end{eqnarray}

Therefore, combining (\ref{pointwise}) and (\ref{spectral}), we have
$$Np^{-\delta N} K \gtrsim \eta^2$$
and so
$$K\gtrsim \eta^2 N^{-1}p^{\delta N}$$
Since this holds for any collection of partition elements of total $\mu$-measure $>\eta$, we conclude that there is at most $\mu$-measure $\eta$ on ergodic components of entropy less than $\delta'\gtrsim \delta>0$.  Taking $\eta\to 0$, we get positive entropy on a.e. ergodic component of $\mu$.

\section*{Acknowledgement}
We thank Nalini Anantharaman, Peter Sarnak, and Akshay Venkatesh for many helpful discussions and encouragement.  Nalini Anantharaman has been very helpful in the writing of the \emph{version fran\c{c}aise abr\'eg\'ee} by correcting and improving a linguistically inadequate original; in preparing this French version we were also assisted by Tony Phillips.

\bibliographystyle{alpha}
\bibliography{my}

\def\cprime{$'$}
\begin{thebibliography}{AKN07}

\bibitem[AKN07]{AKN}
Nalini Anantharaman, Herbert Koch, and St{\'e}phane Nonnenmacher.
\newblock Entropy of eigenfunctions, preprint, 2007.

\bibitem[AN07]{AN}
Nalini Anantharaman and St{\'e}phane Nonnenmacher.
\newblock Half-delocalization of eigenfunctions for the {L}aplacian on an
  {A}nosov manifold.
\newblock {\em Ann. Inst. Fourier (Grenoble)}, 57(7):2465--2523, 2007.
\newblock Festival Yves Colin de Verdi{\`e}re.

\bibitem[Ana08]{Anan}
Nalini Anantharaman.
\newblock Entropy and the localization of eigenfunctions.
\newblock {\em Ann. of Math. (2)}, 168(2):435--475, 2008.

\bibitem[BL03]{BorLin}
Jean Bourgain and Elon Lindenstrauss.
\newblock Entropy of quantum limits.
\newblock {\em Comm. Math. Phys.}, 233(1):153--171, 2003.

\bibitem[BL10]{meElon}
Shimon Brooks and Elon Lindenstrauss.
\newblock Non-localization of eigenfunctions on large regular graphs,
  submitted, 2010.

\bibitem[EKL06]{EKL}
Manfred Einsiedler, Anatole Katok, and Elon Lindenstrauss.
\newblock Invariant measures and the set of exceptions to {L}ittlewood's
  conjecture.
\newblock {\em Ann. of Math. (2)}, 164(2):513--560, 2006.

\bibitem[HS]{Holowinsky-Sound}
Roman Holowinsky and Kannan Soundararajan.
\newblock Mass equidistribution of {H}ecke eigenforms.
\newblock to appear, Annals of Math.

\bibitem[Lin01]{LinHxH}
Elon Lindenstrauss.
\newblock On quantum unique ergodicity for {$\Gamma\backslash\mathbb{H}
  \times\mathbb{H}$}.
\newblock {\em Internat. Math. Res. Notices}, (17):913--933, 2001.

\bibitem[Lin06]{Lin}
Elon Lindenstrauss.
\newblock Invariant measures and arithmetic quantum unique ergodicity.
\newblock {\em Ann. of Math. (2)}, 163(1):165--219, 2006.

\bibitem[Sou09]{Sound}
K.~Soundararajan.
\newblock Quantum unique ergodicity for $sl_2(\mathbb{Z})\backslash\mathbb{H}$,
  preprint, 2009.

\bibitem[SV07]{Lior-Akshay}
Lior Silberman and Akshay Venkatesh.
\newblock On quantum unique ergodicity for locally symmetric spaces.
\newblock {\em Geom. Funct. Anal.}, 17(3):960--998, 2007.

\bibitem[SV10]{Lior-Akshay2}
Lior Silberman and Akshay Venkatesh.
\newblock Entropy bounds for hecke eigenfunctions on division algebras, to
  appear in GAFA, 2010.

\bibitem[Wol01]{Wolpert}
Scott~A. Wolpert.
\newblock Semiclassical limits for the hyperbolic plane.
\newblock {\em Duke Math. J.}, 108(3):449--509, 2001.

\end{thebibliography}

\end{document}